\documentclass[11pt]{article}
\usepackage{amsmath}
\usepackage{amssymb}
\usepackage{graphics}
\usepackage{graphicx}
\def\R{\textbf{R}}
\def\C{\mathbb{C}}%
\def\cp{\mathrm{cap\,}}%
\def\mes{\mathrm{mes\,}}%
\def\m{\mathrm{m\,}}%
\def\vt{\vartheta}%
\def\mathcup{\mathop{\cup}}
\newtheorem{theo}{\hspace*{\parindent}Theorem}
\newtheorem{lemma}{\hspace*{\parindent}Lemma}

\newcounter{theremark}

\title{Two-sided bounds for the logarithmic capacity of multiple intervals}
\author{V.N.\,Dubinin\footnote{Institute of Applied Mathematics, Vladivostok, Russia,
e-mail:\emph{dubinin@iam.dvo.ru}} \,\,and
D.\,Karp\footnote{Institute of Applied Mathematics, Vladivostok,
Russia, e-mail:\emph{dmkrp@yandex.ru}}}
\date{}
\begin{document}
\maketitle
\begin{center}
\parbox{12cm}{
\small\textbf{Abstract.} Potential theory on the complement of a
subset of the real axis attracts a lot of attention both in
function theory and applied sciences. The paper discusses one
aspect of the theory - the logarithmic capacity of closed subsets
of the real line.  We give simple but precise upper and lower
bounds for the logarithmic capacity of multiple intervals and a
lower bound valid also for closed sets comprising an infinite
number of intervals. Using some known methods to compute the exact
values of capacity we demonstrate graphically how our estimates
compare with them. The main machinery behind our results are
separating transformation and dissymmetrization developed by
V.N.\,Dubinin and a version of the latter by K.\,Haliste as well
as some classical symmetrization and projection result for
logarithmic capacity.  The results of the paper improve some
previous achievements by A.Yu.\,Solynin and K.\,Shiefermayr. }
\end{center}

\bigskip


\bigskip


\paragraph{1. Introduction} Let $E$ denote a compact subset of the
complex plane $\C$ and write $g_B(z,\infty)$ for the Green
function of the connected component $B$ of $\C\!\setminus\!{E}$
containing the point at infinity.  \emph{Logarithmic capacity} of
$E$ is defined by
\[
\cp{E}=\exp\{\lim\limits_{z\to\infty}[\log|z|-g_B(z,\infty)]\}.
\]
If $B$ does not admit the Green function we set $\cp{E}=0$.
Logarithmic capacity $\cp{E}$ is equal to the Chebyshev constant
of $E$ and its transfinite diameter \cite{Goluzin,Kirsch}.  Since
logarithmic capacity is not an easy quantity to compute, its lower
and upper estimates are of considerable interest (see, for
instance \cite{Kirsch}).  In this paper we will be concerned with
estimating logarithmic capacity of closed subsets of the real line
in particular those comprising a finite number of intervals. These
type of subsets are obtained, for instance, by Steiner or circular
symmetrization of most one-dimensional sets and hence are extremal
for many problems of function theory. Potential theory on the
complement of such set attracts significant attention
\cite{Andrievskii2,Strang,Totik}. Since $\cp{aE}=|a|\cp{E}$ for
any complex $a$ we may restrict our attention to the subsets of
the interval $[-1,1]$.  The classical bounds in this situation are
\[
(\mes{E})/4\leq\cp{E}\leq{1/2},
\]
where $\mes{E}$ denote the linear Lebesgue measure of $E$.  These
inequalities albeit simple are too rough especially for the sets
consisting of many intervals.  In this connection the question
arises how to obtain more precise estimates taking account of the
structure of $E$ and its dispersion in $[-1,1]$ in terms of
elementary functions.  In the recent work \cite{Schief}
Schiefermayr established the upper bounds for the logarithmic
capacity of $E_{\alpha,\beta}=[-1,\alpha]\cup[\beta,1]$,
$-1<\alpha<\beta<1$,  some of them in terms of elementary
functions.  He also gives a survey of some known and new lower
bounds.  Let us mention some of these bounds together with our
comments and amendments. According to \cite[Theorem~3]{Schief} the
inequality
\begin{equation}\label{eq:ShiefL}
\cp{E_{\alpha,\beta}}\geq
\frac{\sqrt[4]{(1-\alpha^2)(1-\beta^2)}}{\sqrt{(1-\alpha)(1+\beta)}+\sqrt{(1+\alpha)(1-\beta)}}
\end{equation}
holds true with equality attained when $\alpha+\beta=0$.  It is
indicated in \cite{Schief} that Solynin \cite[Section~2.2]{Sol}
proved the lower bound for the logarithmic capacity of multiple
intervals which in the case of two intervals takes the form
\begin{equation}\label{eq:Solynin}
\cp{E_{\alpha,\beta}}\geq
\frac{1}{2}\max\left[\left(\sin\left(\frac{\pi\theta(\beta)}{2\theta(\delta)}\right)\right)^{2\theta(\delta)^2/\pi^2}
\left(\sin\left(\frac{\pi(\pi-\theta(\alpha))}{2(\pi-\theta(\delta))}\right)\right)^{2(\pi-\theta(\delta))^2/\pi^2}\right],
\end{equation}
where here and henceforth $\theta(\gamma)=\arccos(\gamma)$ and the
maximum is taken over all $\delta\in[\alpha,\beta]$.  In view of
the well-known Robinson formula (Lemma~\ref{lm:Rob} below), we
notice that Solynin's inequality is a particular case of the
earlier result of the first author (see, for instance,
\cite[Corollary~1.3]{DubUspehi} and related comments).  Directly
from polarization \cite[Corollary~1.2]{DubUspehi} we get the
following simple upper bound:
\begin{equation}\label{eq:polarization}
\cp{E_{\alpha,\beta}}\leq\cp{E_{-\gamma,\gamma}}=\frac{1}{2}\sqrt{1-\gamma^2}=\frac{1}{4}\sqrt{4-(\alpha-\beta)^2},
\end{equation}
where $\gamma=(\beta-\alpha)/2$.  Another upper bound follows from
an inequality due to Gillis \cite{Gillis}:
\[
\cp{E_{\alpha,\beta}}\leq
2\exp\left[\frac{\log((1+\alpha)/8)\log((1-\beta)/8)}{\log((1+\alpha)(1-\beta)/64)}\right].
\]
The main result of \cite{Schief} is the estimate
\begin{equation}\label{eq:SchiefU}
\cp{E_{\alpha,\beta}}\leq\frac{1+\alpha}{2(1+\beta)}
\exp\left\{2\left(\frac{E}{K}-\frac{(1+\alpha)(1-\beta)}{(1-\alpha)(1+\beta)}\right)\left[\log\frac{\sqrt{2}+\sqrt{1-\alpha}}{\sqrt{1+\alpha}}\right]^2\right\},
\end{equation}
where $K=K(k)$, $E=E(k)$ are Legendre's complete elliptic
integrals of the first and second kinds, respectively,
\[
k=\frac{2(\beta-\alpha)}{(1-\alpha)(1+\beta)},
\]
and it is assumed that $\alpha+\beta\geq{0}$, otherwise one has to
replace $E_{\alpha,\beta}$ with $E_{-\beta,-\alpha}$ having the
same logarithmic capacity. In order to reduce this bound to
elementary functions one may apply two-sided estimates for
elliptic integrals and their ratios from \cite{AQ,AVV}.

In this note we give rather general upper and lower estimates for
the logarithmic capacity of a subset $E$ of $[-1,1]$
(Theorems~\ref{th:1}-\ref{th:upperbound} below).  In particular,
if $E$ consists of $n$ intervals Theorem~\ref{th:Sol-str} gives
the lower bound which coincides with (\ref{eq:Solynin}) for $n=2$
but is stronger than the corresponding result from \cite{Sol} for
$n>2$. This bounds is also stronger than (\ref{eq:ShiefL}) for
$n=2$.  Our upper bound from Theorem~\ref{th:upperbound} is both
very simple and more precise than (\ref{eq:SchiefU}) except for
very narrow neighbourhood of $\beta=1$ where (\ref{eq:SchiefU})
becomes asymptotically precise.  Our bound remains very good for
$n>2$, where it seems to be the only known non-trivial upper
bound.

 The main results of the paper proved
in section~3 are based on the Robinson formula and the estimates
for the capacity of subsets of the unit circle given in section~2.
In the final section~4 we compare our estimates with exact values
computed using the formulas due to Akhieser
\cite{Achieser1,Achieser2} (for $n=2$) and Widom \cite{Widom} (in
a modified form for $n>2$).

\paragraph{2. Auxiliary results.}  Let $\Gamma=\{z:|z|=1\}$.  The following statement
can be derived from the properties of conformal mappings and
symmetry considerations.
\begin{lemma}\emph{(Robinson \cite{Robinson})}.\label{lm:Rob}
Suppose that $F$ is a closed subset of the unit circle $\Gamma$
symmetric with respect to the real axis and let $E$ be its
orthogonal projection onto the real axis.  Then
\[
\cp{E}=\frac{1}{2}(\cp{F})^2.
\]
\end{lemma}

A particular case of the principle of circular symmetrization (see
\cite{DubUspehi}) is the following

\begin{lemma}\emph{(Beurling \cite[p.35-36]{Ahlfors})}.\label{lm:Beu}
The logarithmic capacity of a closed subset of $\Gamma$ having the
length $l$ attains its minimal value $\sin(l/4)$ for a subarc of
$\Gamma$.
\end{lemma}

Define the infinite sectors
$D_k=\{z:\alpha_k<\arg{z}<\alpha_{k+1}\}$, $k=1,2,\ldots,n$,
$\alpha_1<\alpha_2<\cdots<\alpha_n<\alpha_{n+1}=\alpha_1+2\pi$ and
let
$\zeta=p_k(z)=-i(e^{-i\alpha_k}z)^{\pi/(\alpha_{k+1}-\alpha_{k})}$,
$k=1,2,\ldots,n$.  The function $p_k(z)$ effects univalent
conformal mapping of $D_k$ onto the right half-plane
$\Re{\zeta}>0$. For a compact set $F$ satisfying
$F\cap\overline{D_k}\ne\emptyset$, $k=1,2,\ldots,n$, denote by
$F_k$ the union of $p_k(F\cap\overline{D_k})$ with its reflection
with respect to imaginary axis.  According to the terminology of
\cite{DubUspehi}  the family of sets $\{F_k\}_{k=1}^{n}$ is the
result of separating transformation of $F$ with respect to the
family of functions $\{p_k\}_{k=1}^{n}$. A particular case of
\cite[Corollary~1.3]{DubUspehi} is

\begin{lemma}\label{lm:Dub}
The following inequality holds:
\[
\cp{F}\geq\prod\limits_{k=1}^{n}(\cp{F_k})^{(\alpha_{k+1}-\alpha_{k})^2/(2\pi^2)}.
\]
\end{lemma}
Introduce the notation
$F(l,n)=\{z\in\Gamma:|\arg{z}^n|\leq{l/2}\}$, $0<l<2\pi$. The next
statement was proved using dissymmetrization (see
\cite{DubUspehi}):
\begin{lemma}\emph{(Haliste \cite{Haliste})}\label{lm:Hal}
Suppose that $F$ is a union of $n$ closed arcs on the unit circle
$\Gamma$ having total length $l$. Then
\[
\cp{F}\leq\cp{F(l,n)}=[\sin(l/4)]^{1/n}.
\]
\end{lemma}
The last two lemmas also allow for the complete description of
equality cases.

\paragraph{3. Main results}  Given a closed subset $e$ of the interval
$[-1,1]$ put
\[
\m{e}=\int\limits_{e}\frac{dx}{\sqrt{1-x^2}}.
\]
The meaning of this formula in the present context is that
$2\m{e}$ gives the length of the symmetric pre-image of $e$ on
$\Gamma$ under orthogonal projection.  Denote by $E(l,n)$ the
orthogonal projection of $F(l,n)$ onto the real axis.
\begin{theo}\label{th:1}
Let $\{e_k\}_{k=1}^{n}$ be a partitioning of the interval $[-1,1]$
by closed intervals having no common inner points,
$\mathcup\limits_{k=1}^{n}e_k=[-1,1]$.  Then for any closed set
$E\subset[-1,1]$ the inequality
\begin{equation}\label{eq:partition}
\cp{E}\geq\frac{1}{2}\prod\limits_{k=1}^{n}\left[\sin\frac{\pi\m(e_k\cap{E})}{2\m{e_k}}\right]^{2(\m{e_k})^2/\pi^2}
\end{equation}
holds true.  Equality is attained for the sets $E=E(l,n)$,
$0<l<2\pi$, and the partitioning $\{e_k\}_{k=1}^{n}$ by the points
$\cos(\pi{k}/n)$, $k=0,1,\ldots,n$.
\end{theo}
\textbf{Proof.}  First let us consider a partitioning the of unit
circle $\Gamma$ by closed arcs $\sigma_k$ having no common inner
points, $\mathcup\limits_{k=1}^{m}\sigma_k=\Gamma$.  Denote by
$D_k$ the open infinite sector formed by two rays passing through
the endpoints of $\sigma_k$ with the vertex at the origin and
write $\beta_k\pi$ for the angular span of $D_k$ (i.e. Lebesgue
measure of $\sigma_k$).  The function
$p_k(z)=\alpha_kz^{1/\beta_k}$, $|\alpha_k|=1$, effects univalent
conformal mapping of $D_k$ onto the right half plane.  According
to Lemma~\ref{lm:Dub} if a closed subset $F$ of the unit circle
$\Gamma$ intersects all sectors $\overline{D_k}$, then
\[
\cp{F}\geq\prod\limits_{k=1}^{m}(\cp{F_k})^{\beta_k^2/2},
\]
where $\{F_k\}_{k=1}^{m}$ is the result of separating
transformation of $F$ with respect to the family
$\{p_k\}_{k=1}^{m}$.  According to Lemma~\ref{lm:Beu}
\[
\cp{F_k}\geq\sin[(\mes{F_k})/4]=\sin[\mes(\sigma_k\cap{F})/(2\beta_k)],~k=1,\ldots,m.
\]
Hence,
\begin{equation}\label{eq:circle-est}
\cp{F}\geq\prod\limits_{k=1}^{m}\left[\sin\frac{\mes(\sigma_k\cap{F})}{2\beta_k}\right]^{\beta_k^2/2}.
\end{equation}
This estimate reduces to trivially valid inequality
$\cp{F}\geq{0}$ when $F\cap{\overline{D_k}}=\emptyset$ for some
$k$, $1\leq{k}\leq{m}$.

Now let $E$ and $\{e_k\}_{k=1}^{n}$ be the set and the
partitioning from the hypothesis of the theorem.  Let $F$ be the
symmetric subset of the unit circle $\Gamma$ such that it's
orthogonal projection onto the real axis is $E$.  Similarly, the
let $\{\sigma_k\}_{k=1}^{2n}$ be the pre-image of
$\mathcup\limits_{k=1}^{n}e_k$ on $\Gamma$ (the points $+1$ and
$-1$ are assumed to be the endpoint of some arcs $\sigma_k$). Now
as we mentioned earlier if $e_l$ is the projection of $\sigma_k$
then the linear measure of $\sigma_k$ is $\m{e_l}$ and
$\mes(\sigma_k\cap{F})=\m(e_l\cap{E})$.  Inequality
(\ref{eq:circle-est}) applied to the set $F$ and the partitioning
$\{\sigma_k\}_{k=1}^{2n}$ takes the form
\[
\cp{F}\geq\prod\limits_{k=1}^{2n}
\left[\sin\frac{\pi\mes(\sigma_k\cap{F})}{2\mes{\sigma_k}}\right]^{(\mes{\sigma_k})^2/(2\pi^2)}
=\prod\limits_{k=1}^{n}
\left[\sin\frac{\pi\m(e_k\cap{E})}{2\m{e_k}}\right]^{(\m{e_k})^2/\pi^2}.
\]
It is left to apply Lemma~\ref{lm:Rob} which says that
$\cp{E}=(\cp{F})^2/2$.  The case of equality can be verified
directly.~$\square$

Let us notice that the survey paper \cite[page 257]{Kirsch} cites
a particular case of inequality (\ref{eq:circle-est}) with
$\beta_k=2/m$ with incorrectly attributed authorship of this
result (see \cite[Theorem~2.9]{DubUspehi}).

If the set $E$ comprises $n$ intervals and the partitioning
$\{e_k\}_{k=1}^{s}$ is tailored so that $e_k\cap{E}$ consists of
one interval having one common endpoint with $e_k$ then inequality
(\ref{eq:partition}) reproduces inequality (2.7) from \cite{Sol}.
The latter is essentially obtained using the approach similar to
that given in \cite[\S{4}]{DubUspehi} but in terms of reduced
moduli of triangles.  The following theorem gives a strengthening
of this result when $n>2$.

\begin{theo}\label{th:Sol-str}
Suppose $E=\mathcup\limits_{k=1}^{n}[a_k,b_k]$,
$-1=a_1<b_1<a_2<b_2<\cdots<a_n<b_n=1$, $n\geq{2}$.  Then
\begin{equation}\label{eq:Sol-str}
\cp{E}\geq\frac{1}{2}\max\prod\limits_{k=1}^{n}\left\{\frac{1}{2}
\left[\cos\frac{\pi(\theta(b_k)-\theta(\delta_k))}{\theta(\delta_{k-1})-\theta(\delta_k)}
-\cos\frac{\pi(\theta(a_k)-\theta(\delta_k))}{\theta(\delta_{k-1})-\theta(\delta_k)}\right]\right\}^{(\theta(\delta_k)-\theta(\delta_{k-1}))^2/\pi^2},
\end{equation}
where $\delta_0=-1$, $\delta_n=1$ and the maximum is taken over
all $\delta_k$ satisfying $b_k<\delta_k<a_{k+1}$,
$k=1,\ldots,n-1$. The equality is attained for $E=E(l,2(n-1))$,
$0<l<2\pi$, where $E(l,s)$ is defined before Theorem~\ref{th:1},
and the values $\delta_k=\cos(\pi{(n-k)}/n)$, $k=0,1,\ldots,n$.
\end{theo}
\textbf{Proof.}  Let $F$ be the symmetric (with respect to the
real axis) subset of the unit circle $\Gamma$ such that its
orthogonal projection to the real axis is $E$. Let us introduce
the notation $\theta_k=\theta(\delta_{n-k+1})$, $k=1,\ldots,n+1$
and $\theta_k=2\pi-\theta_{2n-k+2}$, $k=n+2,\ldots,2n+1$, where
the numbers $\delta_k$ satisfy the hypotheses of the Theorem.
Applying Lemmas~\ref{lm:Rob} and \ref{lm:Dub} with the values of
$\alpha_k=\theta_k$ and $n$ replaced by $2n$, we will have
\begin{equation}\label{eq:EFk}
\cp{E}=\frac{1}{2}(\cp{F})^2\geq\frac{1}{2}\prod\limits_{k=1}^{2n}(\cp{F_k})^{(\theta_{k+1}-\theta_k)/\pi^2}.
\end{equation}
Each set $F_k$, $k=1,2,\ldots,2n$ comprises one or two subarcs of
$\Gamma$ symmetric with respect to the imaginary axis.  The
orthogonal projection of $F_k$ onto the imaginary axis coincides
with that of $F_{2n-k+1}$ and the length of this projection is
equal to
\[
l_k=\cos\frac{\pi(\theta(b_{n-k+1})-\theta(\delta_{n-k+1}))}{\theta(\delta_{n-k})-\theta(\delta_{n-k+1})}
-\cos\frac{\pi(\theta(a_{n-k+1})-\theta(\delta_{n-k+1}))}{\theta(\delta_{n-k})-\theta(\delta_{n-k+1})},
\]
$k=1,2,\ldots,n$.  In view of Lemma~\ref{lm:Rob} we obtain:
\[
\cp{F_k}=\cp{F_{2n-k+1}}=\sqrt{l_k/2},~~k=1,\ldots,n.
\]
Substituting these values of capacities into (\ref{eq:EFk}) and
changing $n-k+1\mapsto{k}$ we arrive at (\ref{eq:Sol-str}).  The
equality case is straightforward to verify.~~$\square$

As we mentioned earlier for $n=2$  inequalities (\ref{eq:Solynin})
and (\ref{eq:Sol-str}) coincide while for $n\geq{3}$
Theorem~\ref{th:Sol-str} provides more precise estimate than
inequality (2.7) from \cite{Sol} or which amounts to be the same
thing than our inequality (\ref{eq:partition}) with the right
choice of partitioning $\{e_k\}_{k=1}^{n}$. See more details in
section~4 below.

\begin{theo}\label{th:upperbound}
Suppose $E=\mathcup\limits_{k=1}^{n}[a_k,b_k]$,
$-1=a_1<b_1<a_2<b_2<\cdots<a_n<b_n=1$, $n\geq{2}$.  Then
\begin{equation}\label{eq:upperbound}
\cp{E}\leq\frac{1}{2}\left\{\cos\left[\frac{1}{2}\sum\limits_{k=1}^{n-1}(\theta(a_{k+1})-\theta(b_k))\right]\right\}^{1/(n-1)}
\end{equation}
The equality is attained for $E=E(l,2(n-1))$, $0<l<2\pi$.
\end{theo}
\textbf{Proof.} Let the set $F\subset\Gamma$ be the same as in the
proofs of Theorems~\ref{th:1} and \ref{th:Sol-str}.  This set
comprises $2(n-1)$ subarcs of $\Gamma$ with total length
\[
l=2\theta(a_n)+\sum\limits_{k=2}^{n-1}2(\theta(a_k)-\theta(b_k))+2(\pi-\theta(b_1))
=2\pi+2\sum\limits_{k=1}^{n-1}(\theta(a_{k+1})-\theta(b_{k})).
\]
An application of Lemmas~\ref{lm:Rob} and \ref{lm:Hal} yields:
\[
\cp{E}=\frac{1}{2}(\cp{F})^2\leq\frac{1}{2}(\cp{F(l,2(n-1))})^2=\frac{1}{2}(\sin(l/4))^{1/(n-1)}.
\]
The equality case is clear.~$\square$

\paragraph{4. Numeric comparison of estimates.}
In order to compare the presicion of various capacity estimates we
need a method to compute the exact values of the logarithmic
capacity of several intervals.  In case $n=2$ this is provided by
the well known formulas due to Akhieser
\cite{Achieser1,Achieser2}:
\begin{equation}\label{eq:capAchieser}
\cp(E_{\alpha,\beta})=\frac{1}{2}\left[\frac{\vt_4(0;q)\vt_3(0;q)}{\vt_4(\omega;q)\vt_3(\omega;q)}\right]^2,
\end{equation}
where Jacobi's theta functions are
\begin{equation}\label{eq:theta3-def}
\vt_3(z;q)=1+2\sum\limits_{n=1}^{\infty}q^{n^2}\cos(2nz),
\end{equation}
\begin{equation}\label{eq:theta4-def}
\vt_4(z;q)=1+2\sum\limits_{n=1}^{\infty}(-1)^nq^{n^2}\cos(2nz),
\end{equation}
and parameters are found from:
\begin{equation}\label{eq:k}
k^2=\frac{2(\beta-\alpha)}{(1-\alpha)(1+\beta)},~~~k'^2=1-k^2,
\end{equation}
\begin{equation}\label{eq:q-omega}
q=\exp\left(-\pi\frac{K(k')}{K(k)}\right),~~\omega=\frac{\pi
F(\sqrt{(1-\alpha)/2},k)}{2{K(k)}},
\end{equation}
where $F(\lambda,k)$ is the first incomplete elliptic integral of
Legendre.

Formula (\ref{eq:capAchieser}) was generalized to three intervals
by Falliero and Sebbar \cite{Falliero,Falliero2} in terms of genus
2 theta functions.  For arbitrary $n$ a mehod to compute the
capacity via Schwarz-Christoffel map was first given by Widom in
\cite{Widom}.  We will use a slightly different guise of his
formula (see further development of Widom's idea in
\cite{Embree}). Indeed, one can verify directly that the Green
function of $\C\backslash{E}$ with pole at infinity is (see
\cite{Andrievskii2,Embree,Strang,Widom}):
\begin{equation}\label{eq:Green-multiple}
g(z)=\Re{F(z)},~~~F(z)=\int\limits_{a_1}^z\frac{p(t)dt}{\sqrt{q(t)}},
\end{equation}
where
\[
q(t)=\prod\limits_{i=1}^n(t-a_i)(t-b_i)
\]
and the branch of square root $\sqrt{q(t)}$ is chosen so that it
is asymptotically $t^n$ near infinity.  The polynomial
\[
p(t)=t^{n-1}+c_{n-2}t^{n-2}+\cdots+c_0
\]
is chosen so that the Schwarz-Christoffel map $F(z)$ maps
$(a_i,b_i)$ into imaginary axis. Since $q(t)>0$ on
$\R\backslash{E}$ we have the linear system of equations
$F(b_i)=F(a_{i+1})$, $i=1,\ldots,n-1$, or
\[
\int\limits_{b_i}^{a_{i+1}}\frac{p(t)dt}{\sqrt{q(t)}}=0,~~~
i=1,\ldots,n-1,
\]
for the definition of the coefficients $c_k$, $k=0,1,\ldots,n-2$.
Further, the Green function has the expansion
\[
g(z)=\ln|z|+R+o(1),~~~~|z|\to\infty,
\]
where $R$ is the Robin constant of $E$ and $\cp{E}=e^{-R}$. Since
the asymptotic expansion of the Green function is true no matter
from which direction we approach infinity, we can move along the
real axis:
\[
R=\lim\limits_{x\to\infty}\left(\int\limits_{b_n}^x\frac{p(t)dt}{\sqrt{q(t)}}-\ln(x)\right)
=\int\limits_{b_n}^{\infty}\left[\frac{p(t)}{\sqrt{q(t)}}-\frac{1}{t}\right]dt=
\int\limits_{b_n}^{\infty}\frac{tp(t)-\sqrt{q(t)}}{t\sqrt{q(t)}}dt,
\]
since
\[
\Re\left[\int\limits_{a_1}^{b_n}\frac{p(t)dt}{\sqrt{q(t)}}\right]=0.
\]
So, finally
\begin{equation}\label{eq:capE}
\cp{E}=\exp\left(\int\limits_{b_n}^{\infty}\left[\frac{1}{t}-\frac{p(t)}{\sqrt{q(t)}}\right]dt\right)
=\exp\left(\int\limits_{b_n}^{\infty}\frac{\sqrt{q(t)}-tp(t)}{t\sqrt{q(t)}}dt\right).
\end{equation}

\begin{figure}[ht]
\hspace*{0.2in}\begin{tabular}{ccc}
\includegraphics[angle=0,width=80mm]{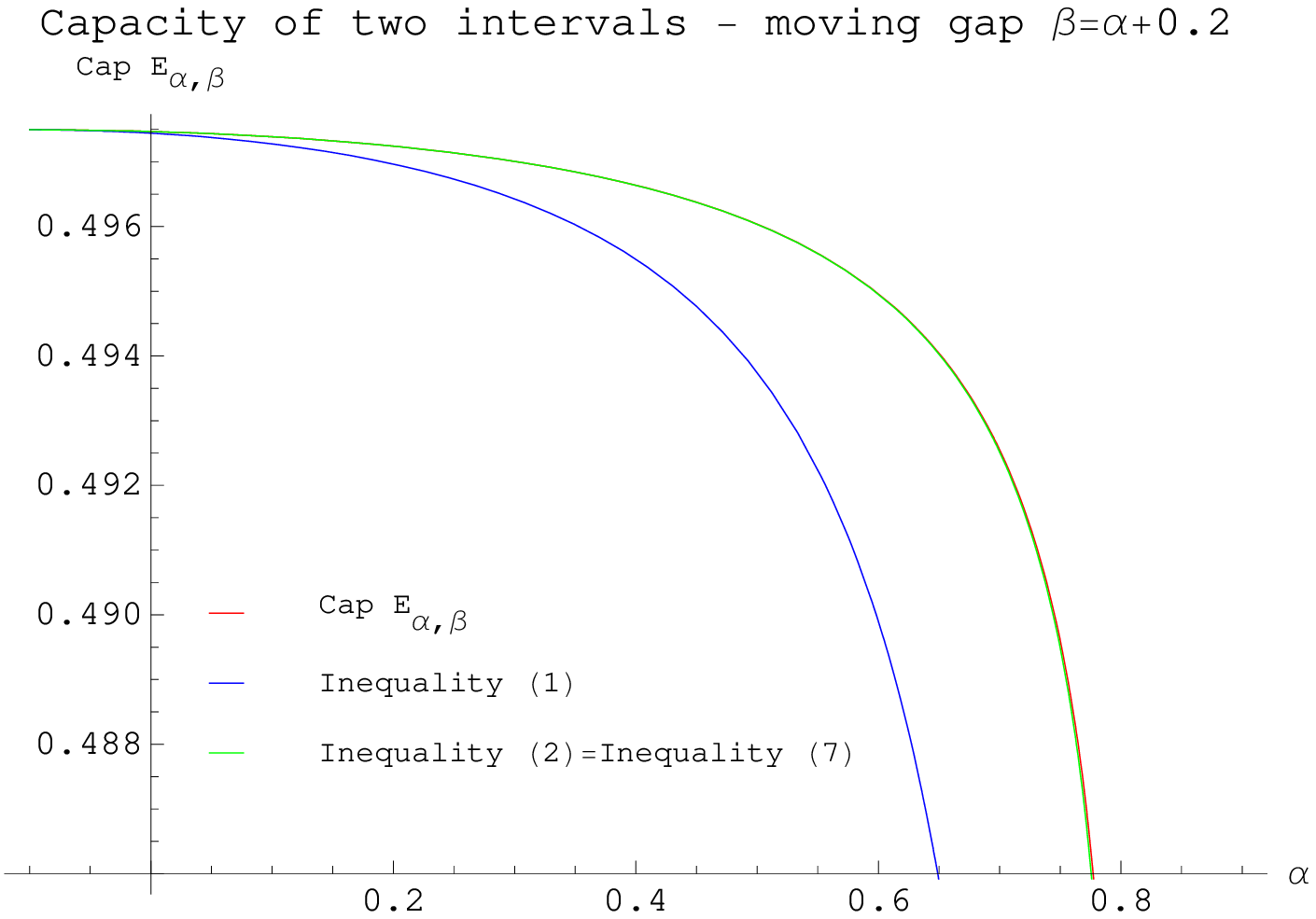}&\hspace*{-0.3in}
\includegraphics[angle=0,width=80mm]{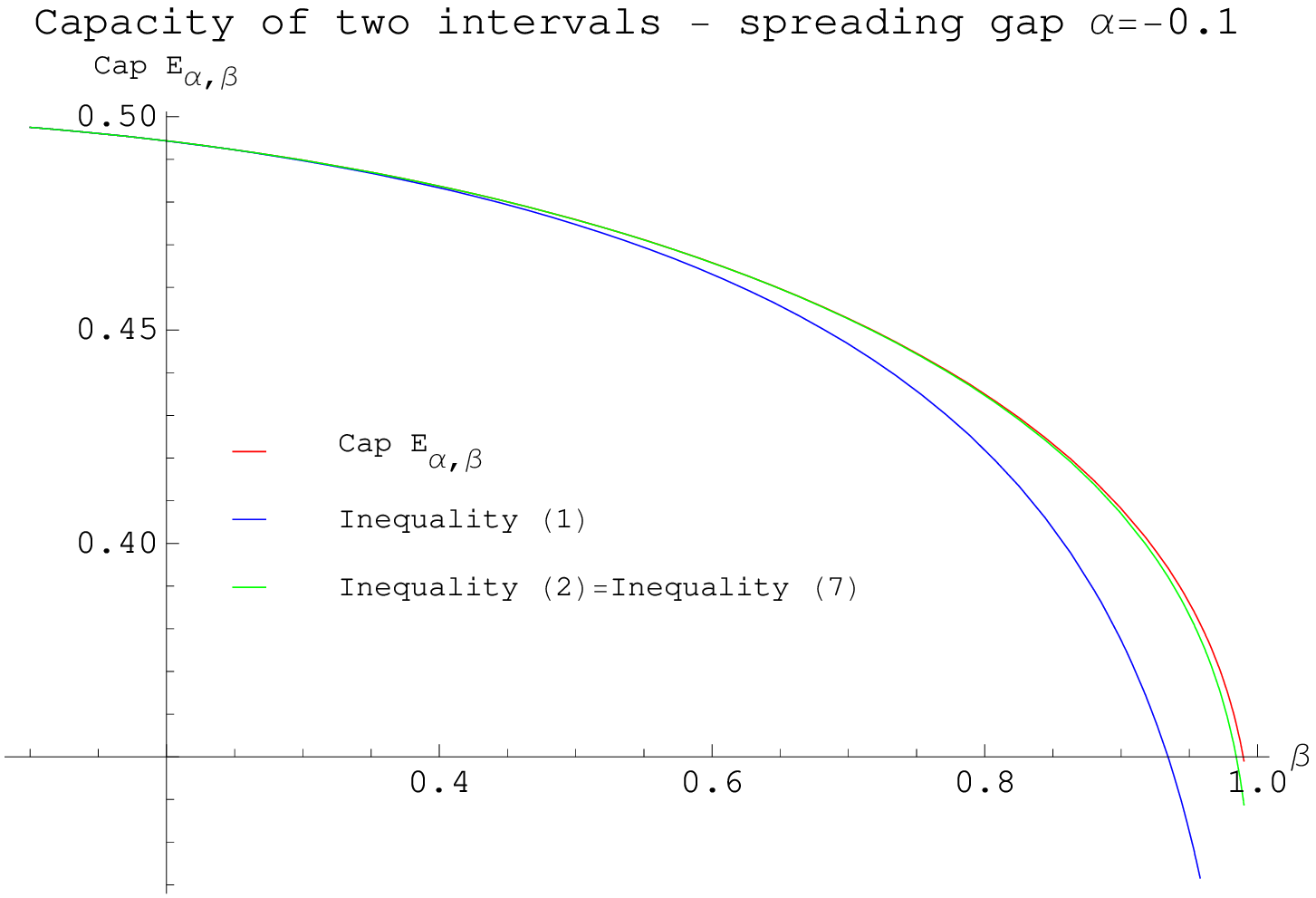}&\hspace*{-0.3in}
\end{tabular}
\caption{The capacity of two intervals and its lower estimates in
two situations}\label{fig:lower2}
\end{figure}

We will make some numerical comparisons between various estimates
and precise capacity values computed using formula
(\ref{eq:capAchieser}) for $n=2$ and formula (\ref{eq:capE}) for
$n=3$.   For convenience let us record the lower bound from
\cite{Sol} and our bounds from Theorems~\ref{th:Sol-str} and
\ref{th:upperbound} for these values of $n$.  For $n=2$, we have
as before $E_{\alpha,\beta}=[-1,\alpha]\cup[\beta,1]$,
$-1<\alpha<\beta<1$, and the lower bounds due to Schiefermayr
\cite{Schief} and Solynin \cite{Sol} are given by formulas
(\ref{eq:ShiefL}) and (\ref{eq:Solynin}), respectively.  One can
verify by straightforward computation that  for $n=2$ our lower
bound (\ref{eq:Sol-str}) reduces to (\ref{eq:Solynin}).  We
demonstrate these bounds in two typical situations - that of a
moving gap and that of a spreading gap on Figure~\ref{fig:lower2}.
Clearly in both situation Solynin's bound (\ref{eq:Solynin}) or
equivalently our bound (\ref{eq:Sol-str}) provide a better
estimate.  Note also that for the moving gap situation the
capacity is monotone and maximal when $\alpha=-\beta$. A proof of
this and similar facts and their generalizations can be found in
\cite{DubKarp}.

\begin{figure}[ht]
\hspace*{0.2in}\begin{tabular}{ccc}
\includegraphics[angle=0,width=80mm]{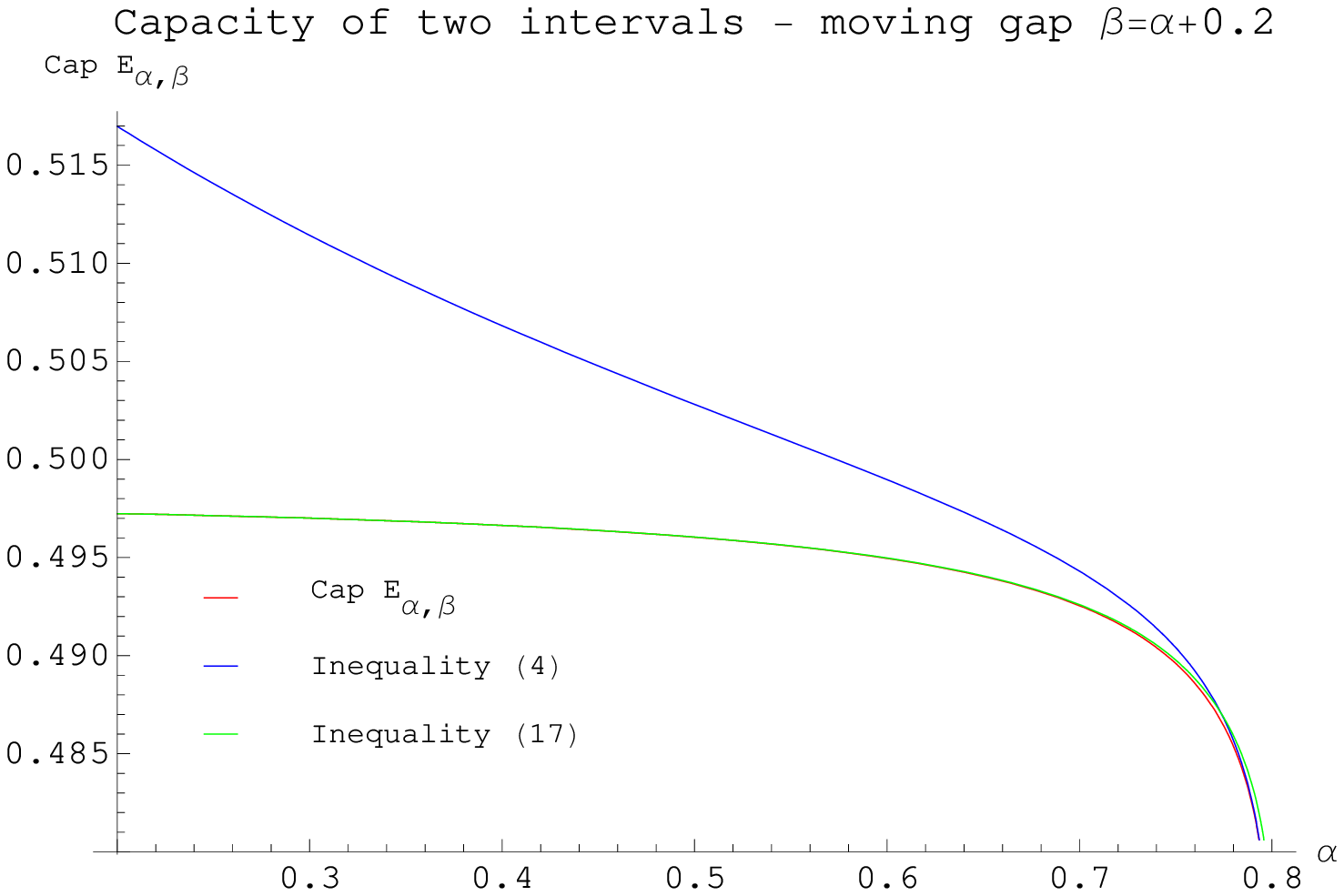}&\hspace*{-0.3in}
\includegraphics[angle=0,width=80mm]{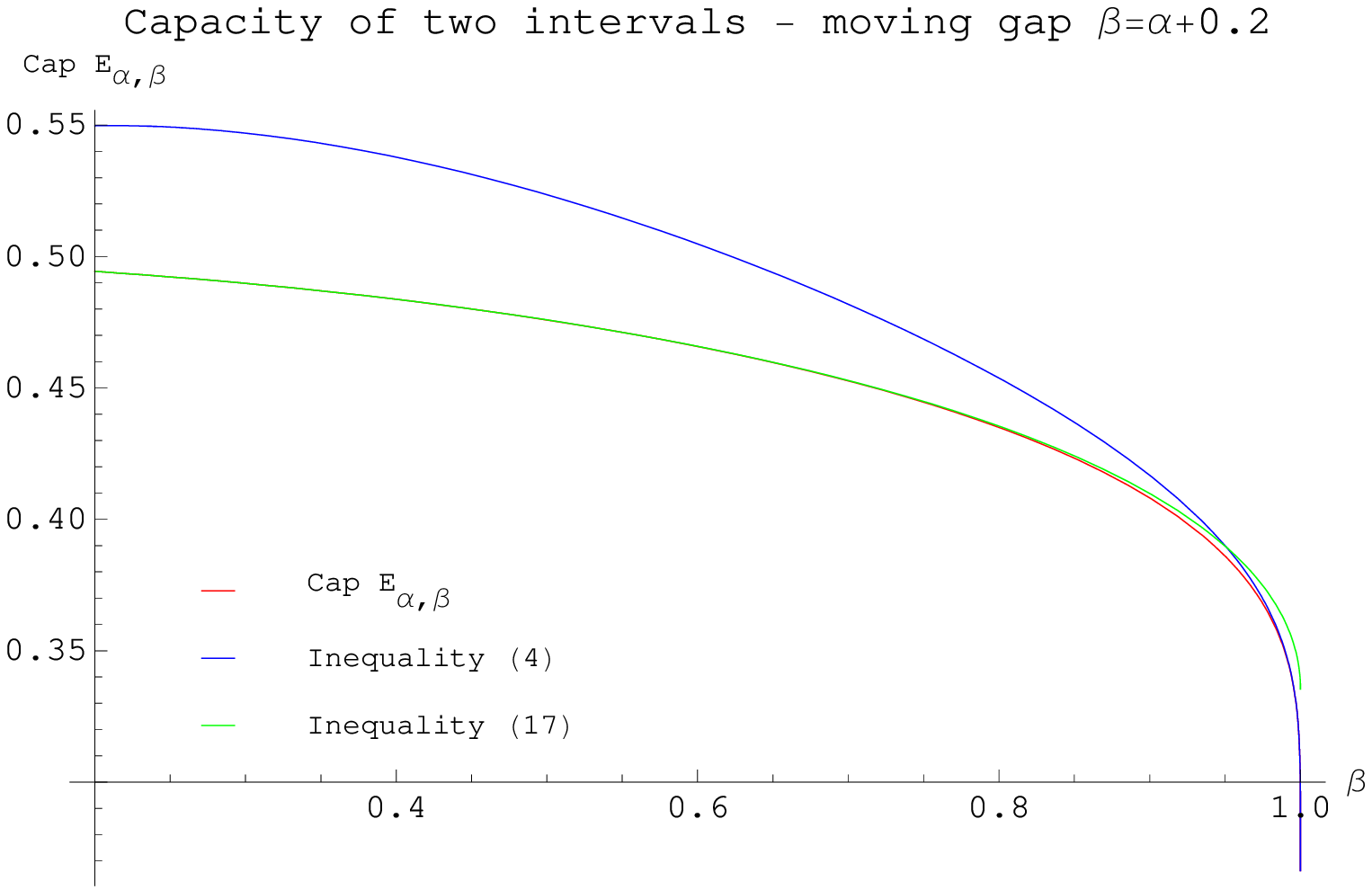}&\hspace*{-0.3in}
\end{tabular}
\caption{The capacity of two intervals and its upper estimates in
two situations}\label{fig:upper2}
\end{figure}

The upper bound of Schiefermayr \cite{Schief} and our simple upper
bound obtained by polarization are given by (\ref{eq:SchiefU}) and
(\ref{eq:polarization}).  The upper bound (\ref{eq:upperbound})
for $n=2$ reads (recall that
$\theta(\alpha)\equiv\arccos(\alpha)$):
\begin{equation}\label{eq:twoU}
\cp{E}\leq\frac{1}{2}\cos\left([\theta(\beta)-\theta(\alpha)]/2\right).
\end{equation}
Figure~\ref{fig:upper2} illustrates these two bounds and the value
of capacity computed by (\ref{eq:capAchieser}). Again we have
chosen two typical situation - the moving gap and the spreading
gap.  Since Schiefermayr's inequality (\ref{eq:SchiefU}) is
asympotically precise when one of the intervals vanishes his bound
becomes more precise in a narrow neighbourhoods of the right end
of the $\alpha$ and $\beta$ ranges.  Our bound provides almost
uniform and tight fit in the whole parameter ranges.

\bigskip

\bigskip

\begin{figure}[ht]
\hspace*{0.2in}\begin{tabular}{ccc}
\includegraphics[angle=0,width=80mm]{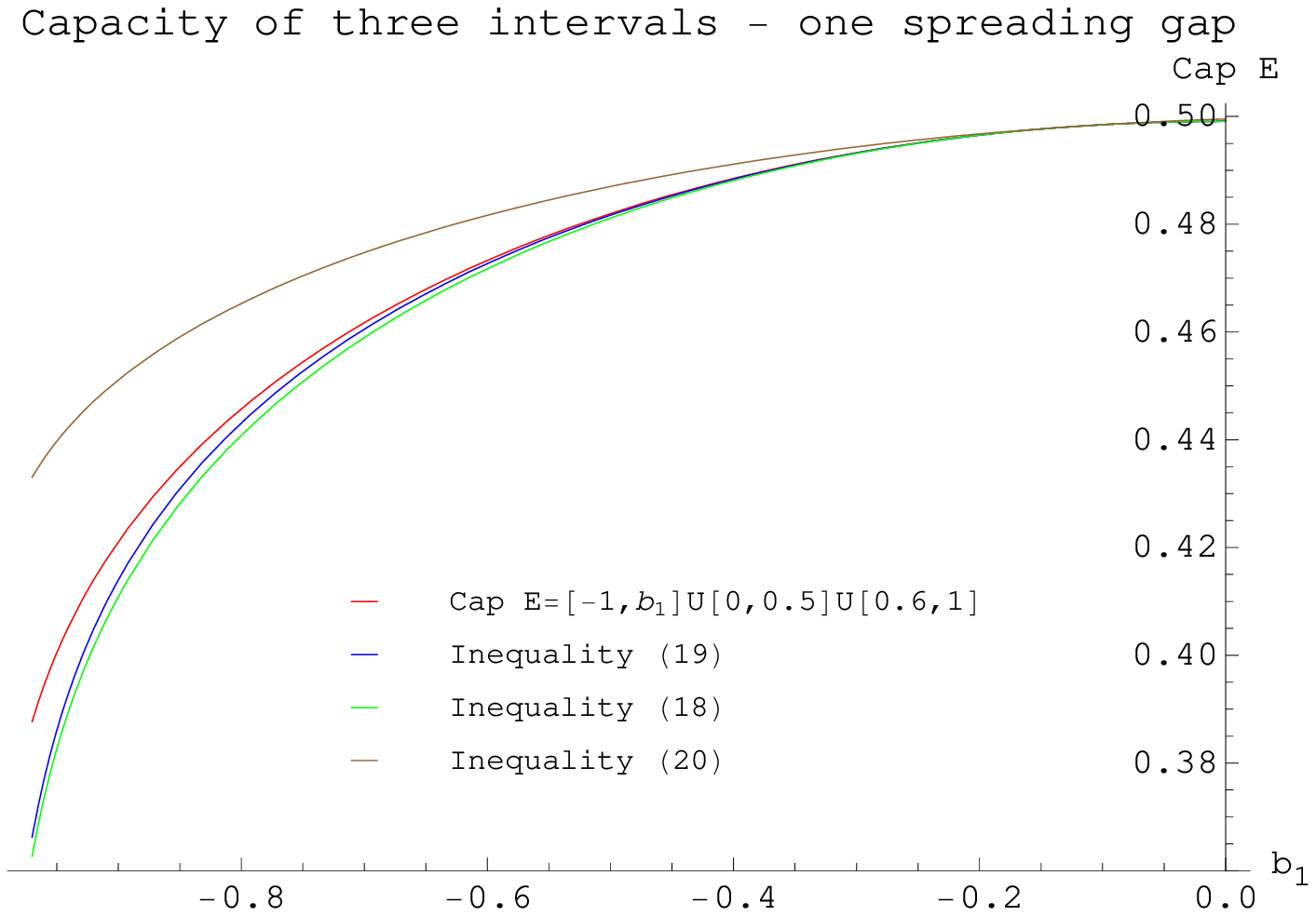}&\hspace*{-0.3in}
\includegraphics[angle=0,width=80mm]{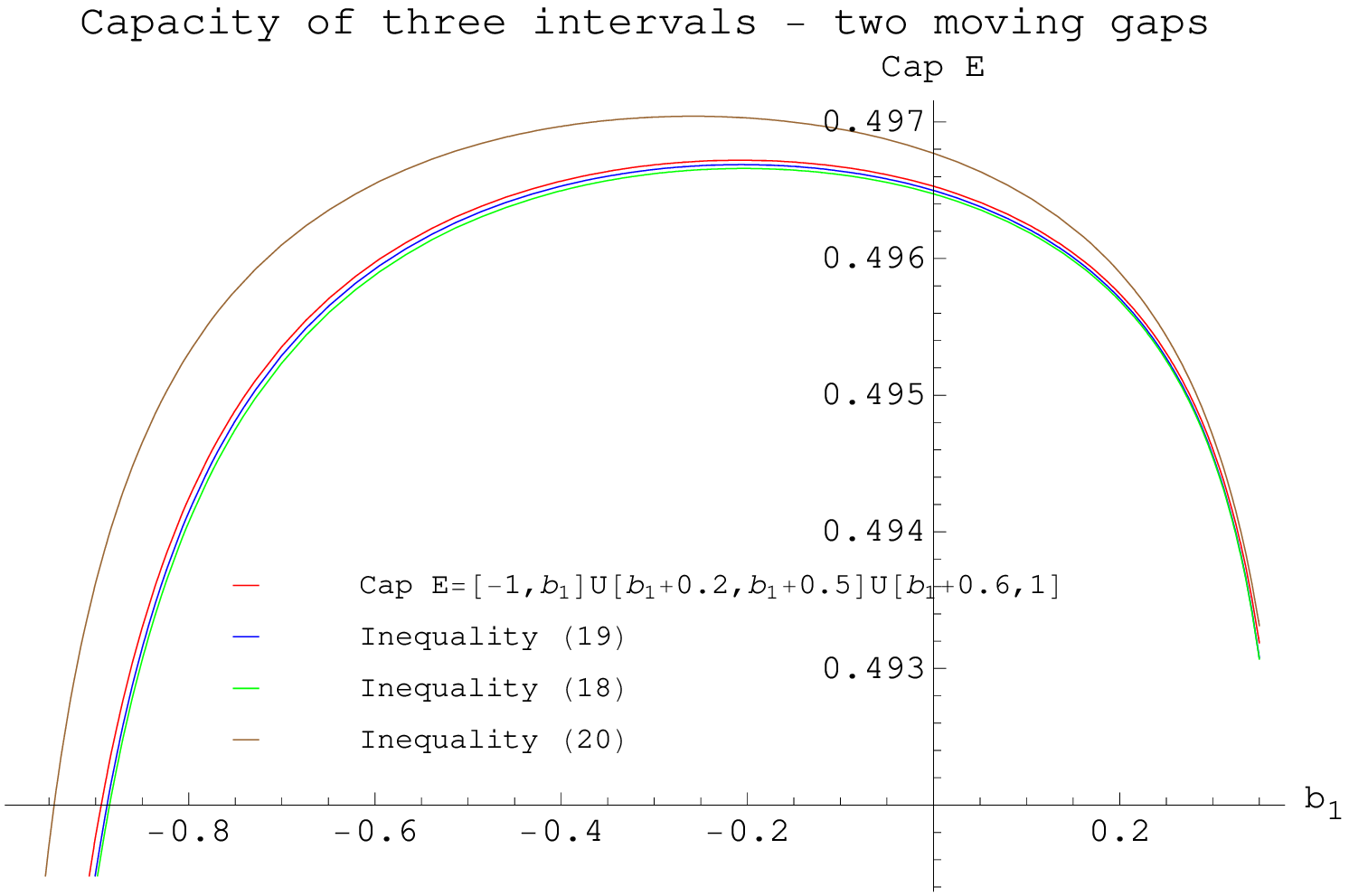}&\hspace*{-0.3in}
\end{tabular}
\caption{The capacity of three intervals and its lower and upper
estimates in two situations}\label{fig:3int}
\end{figure}

For $n>2$ our lower bound (\ref{eq:Sol-str}) differs from that of
Solynin \cite[formula (2.7)]{Sol}. In particular for $n=3$
\cite[formula (2.7)]{Sol} can be written as
\begin{multline}\label{eq:Sol3int}
\cp{E}\geq\frac{1}{2}\max\left(\sin\frac{\pi(\pi-\theta(b_1))}{2(\pi-\theta(\delta_1))}\right)^{2(\pi-\theta(\delta_1))^2/\pi^2}
\left(\sin\frac{\pi(\theta(a_2)-\theta(\gamma))}{2(\theta(\delta_1)-\theta(\gamma))}\right)^{2(\theta(\delta_1)-\theta(\gamma))^2/\pi^2}
\\
\times\left(\sin\frac{\pi(\theta(\gamma)-\theta(b_2))}{2(\theta(\gamma)-\theta(\delta_2))}\right)^{2(\theta(\gamma)-\theta(\delta_2))^2/\pi^2}
\left(\sin\frac{\pi\theta(a_3)}{2\theta(\delta_2)}\right)^{2\theta(\delta_2)^2/\pi^2},
\end{multline}
where the maximum is taken over $\delta_1\in(b_1,a_2)$,
$\gamma\in(a_2,b_2)$ and $\delta_2\in(b_2,a_3)$, while the bound
(\ref{eq:Sol-str}) takes the form
\begin{multline}\label{eq:Dub3int}
\cp{E}\geq\frac{1}{2}\max\left(\frac{1}{2}
\left[1+\cos\frac{\pi(\theta(b_1)-\theta(\delta_1))}{\pi-\theta(\delta_1)}\right]\right)^{(\pi-\theta(\delta_1))^2/\pi^2}
\\
\times\left(\frac{1}{2}
\left[\cos\frac{\pi(\theta(b_2)-\theta(\delta_2))}{\theta(\delta_{1})-\theta(\delta_2)}
-\cos\frac{\pi(\theta(a_2)-\theta(\delta_2))}{\theta(\delta_{1})-\theta(\delta_2)}\right]\right)^{(\theta(\delta_{1})-\theta(\delta_2))^2/\pi^2}
\left(\frac{1}{2}
\left[1-\cos\frac{\pi\theta(a_3)}{\theta(\delta_{2})}\right]\right)^{\theta(\delta_{2})^2/\pi^2},
\end{multline}
where the maximum is taken over $\delta_1\in(b_1,a_2)$ and
$\delta_2\in(b_2,a_3)$. It can be seen from the proof of
\cite[formula (2.7)]{Sol} that for any choice of $\gamma$ the
latter estimate is greater (i.e. better) than the former if
$\delta_i$, $i=1,2$, take the same values in both formulas.
Similar statement holds for $n>3$.

The upper bound (\ref{eq:upperbound}) for $n=3$ reads
\begin{equation}\label{eq:uperr3int}
\cp{E}\leq\frac{1}{2}\left\{\cos\left([\theta(a_{2})-\theta(b_1)+\theta(a_{3})-\theta(b_2)]/2\right)\right\}^{1/2}.
\end{equation}

We compare all three bounds with each other and the value of
capacity computed by (\ref{eq:capE}) on Figure~\ref{fig:3int}.
Again we have chosen a spreading gap as one typical situation
while as the other typical situation we have taken two
simultaneously moving gaps. The figure confirms our prediction
that that bound (\ref{eq:Dub3int}) is more precise than
(\ref{eq:Sol3int}) in both situation.

\paragraph{5. Acknowledgements.}  This work supported by Far Eastern Branch of the Russian Academy of
Sciences (grants 09-III-A-01-008 and 09-II-CO-01-003), Russian
Basic Research Fund (grant 08-01-00028-a) and the Presidential
Grant for Leading Scientific Schools (grant 2810.2008.1).


\begin{thebibliography}{99}
\bibitem{Achieser1} N.I.\,Achieser, \"{U}ber einige Funktionen, welche in zwei
gegebenen Intervallen am weginsten von Null abweichen, I. Teil,
\emph{Bulletin de Academie des Sciences de L'URSS}, 1932,
1163--1202.
\bibitem{Achieser2} N.I.\,Achieser, \"{U}ber einige Funktionen, welche in zwei
gegebenen Intervallen am weginsten von Null abweichen, II. Teil,
\emph{Bulletin de Academie des Sciences de L'URSS}, 1933,
309--344.
\bibitem{AQ} H.\,Alzer and S.-L.\,Qui, Monotonicity theorems and
inequalities for complete elliptic integrals, \emph{Journal of
Comp. and Appl. Math.}, 172, 2004, 289--312.
\bibitem{AVV} G.D.\,Anderson, M.K.\,Vamanamurthy and
M.\,Vourinen, Functional inequalities for complete elliptic
integrals and their ratios, \emph{SIAM J. Math. Anal.}
\textbf{21}, no.2 (1990), 536--549.
\bibitem{Andrievskii2} V.V.\,Andrievskii,  On the Green function for a complement of a finite number of real intervals.
\emph{Constr. Approx.}, \textbf{20} (2004), 4, 565--583.
\bibitem{Ahlfors} L.V.\,Ahlfors, \emph{Conformal invariants. Topics in geometric function theory}, New York,
McGrow-Hill Book Co., 1973.
\bibitem{DubUspehi}V.N.\,Dubinin, Symmetrization in the geometric theory of functions of a complex variable.
\emph{Russ. Math. Surv.} vol.\textbf{49}, no.1, 1--79 (1994);
translation from \emph{Usp. Mat. Nauk} \textbf{49}, No.1(295),
3--76 (1994).
\bibitem{DubKarp} V.N.\,Dubinin and  D.\,Karp,
Capacities of certain plane condensers and sets under simple
geometric transformations,  \emph{Complex Variables and Elliptic
Equations}, vol.\textbf{53}, no. 6 (2008), 607--622,
http://dx.doi.org/10.1080/17476930701734292
\bibitem{Embree} M.\,Ebree, L.N.\,Trefethen, Green's functions for
multiply connected domains via conformal mapping, \emph{SIAM
Review}, vol.\textbf{41}, no. 4 (1999), 745--761.
\bibitem{Falliero}T. Falliero, A. Sebbar, Capacit\'{e} de la
r\'{e}union de troi intervalles et function th$\mathrm{\hat{e}}$ta
de genre 2, \emph{C.R. Acad.Sci. Paris}, \textbf{328}, S\'{e}rie
I, 763--766, 1999.
\bibitem{Falliero2} T.\,Falliero and A.\,Sebbar, Capacite de la reunion de trois
intervalles et fonctions theta de genre 2. \emph{J. Math. Pures
Appl.} (9) 80 (2001),  4, 409--443.
\bibitem{Gillis} J.\,Gillis, Tchebycheff polynomials and the
transifinite diameter, \emph{American Journal of Mathematics}, 63,
1941, 283--290.
\bibitem{Goluzin} G.M.\,Goluzin,
\emph{Geometric theory of functions of a complex variable},
Providence, R. I.:American Mathematical Society (AMS), VI, 1969.
\bibitem{Haliste}K.\,Haliste, On extremal configuration for
capacity, \emph{Arkiv for Mat.}, vol.\textbf{27}, no.1, 1989,
97--104.
\bibitem{Kirsch}S.\,Kirsch, Transfinite diameter, Chebyshev Constant and
Capacity, \emph{Handbook of Complex Ananlysis: Geometric Function
Theory, Volume~2}, ed. by R.\,K\"{u}hnau, Elsevier, 2005.
\bibitem{Robinson} R.M.\,Robinson, On the transifinite diameters
of some related sets, \emph{Mathematische Zeitschrift}, 108, 1969,
377--380.
\bibitem{Schief} K.\,Schiefermayr, An upper bound for the
logarithmic capacity of two intervals, \emph{Complex Variables and
Elliptic Equations}, vol.\textbf{53}, no.1(2008), 65--75.
\bibitem{Strang} J.\,Shen, G.\,Strang and A.J.\,Wathen,  The potential
theory of several intervals and its applications. \emph{Appl.
Math. Optim.} vol.\textbf{44}, no. 1(2001), 67--85.
\bibitem{Sol} A.Yu.\,Solynin, Extremal configurations in some
problems of capacity and harmonic measure, \emph{Journal of
Mathematical Sciences}, 89, 1998, 1031--1049.
\bibitem{Totik} V.\,Totik, \emph{Metric properties of harmonic measure},
Memoirs of AMS, vol.\textbf{184} (2006), no.867.
\bibitem{Widom} H.\,Widom, Extremal polynomials associated with a
system of curves in the complex plane, \emph{Adv. in Math.},
\textbf{3}, 127--232.



\end{thebibliography}
\end{document}